\pdfoutput=1
\RequirePackage{ifpdf}
\ifpdf 
\documentclass[pdftex]{sigma}
\else
\documentclass{sigma}
\fi

\usepackage{enumerate}
\usepackage[all]{xy}
\usepackage{pb-diagram,pb-xy}

\begin{document}

\renewcommand{\PaperNumber}{006}

\FirstPageHeading

\ShortArticleName{On Lie Algebroids and Poisson Algebras}

\ArticleName{On Lie Algebroids and Poisson Algebras}

\Author{Dennise GARC\'IA-BELTR\'AN~$^\dag$, Jos\'e A. VALLEJO~$^\dag$ and Yuri\u{\i} VOROBJEV~$^\ddag$}

\AuthorNameForHeading{D.~Garc\'ia-Beltr\'an, J.A.~Vallejo and Yu.~Vorobjev}

\Address{$^\dag$~Facultad de Ciencias, Universidad Aut\'onoma de San Luis Potos\'i, M\'exico}
\EmailD{\href{mailto:dennise.gb@alumnos.uaslp.edu.mx}{dennise.gb@alumnos.uaslp.edu.mx}, \href{mailto:jvallejo@fc.uaslp.mx}{jvallejo@fc.uaslp.mx}}
\URLaddressD{\url{http://galia.fc.uaslp.mx/~jvallejo/}}

\Address{$^\ddag$~Departamento de Matem\'aticas, Universidad de Sonora, M\'exico}
\EmailD{\href{mailto:yurimv@guaymas.uson.mx}{yurimv@guaymas.uson.mx}}

\ArticleDates{Received June 08, 2011, in f\/inal form February 01, 2012; Published online February 10, 2012}

\Abstract{We introduce and study a class of Lie algebroids associated to
faithful modules which is motivated by the notion of cotangent Lie algebroids of
Poisson manifolds. We also give a classif\/ication of transitive Lie algebroids and
describe Poisson algebras by using the notions of algebroid and Lie connections.}

\Keywords{transitive Lie algebroids; Lie--Rinehart algebras; Poisson brackets; algebraic con\-nec\-tions}

\Classification{53D17; 17B63}

\section{Introduction}

Let $(M,\{ \cdot,\cdot\})$ be a Poisson manifold equipped with the bracket $\{ \cdot,\cdot\}$, which is determined by a~Poisson bivector $\mathbf{P}\in \Gamma \Lambda^2 TM$. It is well known that the cotangent bundle $T^* M$ carries a~natural Lie algebroid structure (see Section \ref{Preliminaries} for def\/initions), that is, on dif\/ferential $1$-forms (sections of~$T^*M$) the following bracket can be def\/ined (see~\cite[Proposition~14.19]{AM69}),
\begin{gather}\label{eq1}
[\![ \alpha ,\beta]\!] =i_{P\alpha}\mathrm{d}\beta -i_{P\beta}\mathrm{d}\alpha +\mathrm{d}(\beta (P\alpha)),
\end{gather}
where $P:T^* M\rightarrow TM$ is the vector bundle morphism canonically induced by $\mathbf{P}$ through $\beta (P\alpha )=\mathbf{P}(\alpha ,\beta )$. Moreover, the following Leibniz-like rule is satisf\/ied, for any smooth function $f\in \mathcal{C}^{\infty}(M)$:
\begin{gather*}
[\![ \alpha ,f\beta ]\!] = f[\![ \alpha ,\beta]\!] +(P\alpha )(f)\cdot \beta .
\end{gather*}
That means that the corresponding anchor map is just $P$.
Note that in the case when~$\alpha$,~$\beta$ in~\eqref{eq1} are exact ($\alpha =\mathrm{d}f$ and $\beta =\mathrm{d}g$ for some $f,g\in \mathcal{C}^{\infty}(M)$), we have $\mathbf{P}(\alpha ,\beta )=\{ f,g\}$ and formula~\eqref{eq1} reads:
\begin{gather*}
[\![ \mathrm{d}f,\mathrm{d}g ]\!] =\mathrm{d}\{ f,g\}.
\end{gather*}
Moreover, the bracket \eqref{eq1} of two closed forms is again closed.

In fact, the properties above characterize the Lie algebroid structure on the cotangent bundle that comes from a Poisson bracket on the base manifold (see \cite{Hue90,Du05}). More precisely, given a~Lie algebroid $(T^* M, [\![ \cdot ,\cdot ]\!] ,\rho )$ on the cotangent bundle $T^* M$, there exists a Poisson bracket~$\{ \cdot,\cdot\}$ on~$M$ such that $\rho =P$ if and only if the following conditions are satisf\/ied:
\begin{enumerate}[(a)]\itemsep=0pt
\item\label{a} $\rho$ is skew-symmetric, i.e.\ $\beta (\rho (\alpha ))=-\alpha (\rho (\beta ))$, for all $\alpha,\beta\in\Omega^1 (M)$.
\item\label{b} If $\alpha,\beta\in\Omega^1 (M)$ are closed, then $[\![ \alpha ,\beta]\!]$ is also closed.
\end{enumerate}
The aim of the present work is to put the study of cotangent Lie algebroids of Poisson manifolds in an algebraic framework.
In the f\/irst part (Section~\ref{LiePoisson}), we study and characterize a class of Lie algebroids with properties similar to~\eqref{a},~\eqref{b}, which we call Lie algebroids of Poisson type. To this end, we will work in a slightly more general context than that of vector bundles over manifolds, using the notion of Lie--Rinehart algebra (see~\cite{Hue90,Hue98}), also called a Lie pseudoalgebra (see \cite{KM90, Mc95, Mc05, Ri63} and, particularly, \cite{KM02} for many interesting remarks on the evolution of these notions and as a general reference), although we will continue using the denomination ``Lie algebroid'' for them (as they are the algebraic version of the geometric Lie algebroids). The dictionary
in Example \ref{ec} can be used at any time to obtain the corresponding expressions for geometric Lie algebroids. We give several examples illustrating the dif\/ferent situations that can appear.

In the second part (Sections \ref{Transitive} and \ref{Poisson}), we deepen in the relationship between transitive algebroids and Poisson structures for a certain class of spaces, those of the form \mbox{$\mathrm{Der}(\mathcal{A})\oplus V$}, where $\mathcal{A}$ is a commutative algebra and $V$ an $\mathcal{A}$-module. We describe parametrizations of the transitive Lie algebroids on $\mathrm{Der}(\mathcal{A})\oplus V$ following the techniques exposed in \cite{Yu01, Yu04}, which are based on the use of a connection on a Lie algebroid. For completeness, we include a subsection in the preliminaries devoted to the topic of connections in an algebraic setting. Once the parametrization is given, we apply it to prove that a transitive algebroid endowed with a connection is isomorphic to one of the form $\mathrm{Der}(\mathcal{A})\oplus V$. Finally, we obtain new classes of Poisson algebras on $\mathcal{A}\oplus V$ starting from Poisson algebras on $\mathcal{A}$.

\section{Preliminaries}\label{Preliminaries}

Throughout the document, unless otherwise explicitly stated, $\mathcal{A}$ denotes an associative, commutative algebra with identity element $1_{\mathcal{A}}$, over a commutative ring $\mathcal{R}$ with identity element~$1_{\mathcal{R}}$.

\subsection{Derivations and connections in commutative algebras}\label{subsec21}
In subsequent sections, we will need to introduce connections on an algebroid. In our algebraic setting, the most appropriate notion of connection is Koszul's one, which is given in terms of derivations.
\begin{definition}
A derivation of the algebra $\mathcal{A}$ over $\mathcal{R}$ is a map $X \in \mathrm{Hom}_{\mathcal{R}}(\mathcal{A},\mathcal{A})$ satisfying the Leibniz rule
\begin{gather*}
X (f\cdot g)=X(f)\cdot g+f\cdot X(g).
\end{gather*}
The set of derivations of $\mathcal{A}$ over $\mathcal{R}$ is denoted $\mathrm{Der}_{\mathcal{R}}(\mathcal{A})$ or simply
$\mathrm{Der}(\mathcal{A})$ when there is no risk of confusion about the ring $\mathcal{R}$.
\end{definition}
\begin{remark}
The def\/inition just given can be extended to the case of derivations of $\mathcal{A}$ over $\mathcal{R}$ with values in an $\mathcal{A}$-module $M$. These are Abelian groups morphisms $X:\mathcal{A}\rightarrow M$ satisfying the Leibniz rule above, and form an $\mathcal{A}$-module denoted $\mathrm{Der}_{\mathcal{R}}(\mathcal{A},M)$ or simply $\mathrm{Der}(\mathcal{A},M)$.
\end{remark}

 The set $\mathrm{Der}(\mathcal{A})$ has an $\mathcal{R}$-Lie algebra structure when endowed with the commutator of endomorphisms, given by $[X ,Y ]=X \circ Y -Y \circ X$.

Note also that, if $\mathcal{A}$ as an $\mathcal{R}$-module is faithful, then for every $X\in \mathrm{Der}(\mathcal{A})$ we have $X (1_{\mathcal{A}})=0$ and indeed $X(r)=0$ for every $r\in \mathcal{R}$ viewed as a subalgebra of $\mathcal{A}$.
\begin{definition}
Let $M$ be an unitary $\mathcal{A}$-module. A derivation law, or Koszul connection, on $M$ is an $\mathcal{A}$-linear mapping $\nabla :\mathrm{Der}(\mathcal{A})\rightarrow \mathrm{Hom}_{\mathcal{R}}(M,M)$ (the image of $X\in \mathrm{Der}(\mathcal{A})$ denoted $\nabla_X$) such that
\begin{gather*}
\nabla_X (f\cdot m)=X(f)\cdot m+f\cdot \nabla_X(m).
\end{gather*}
\end{definition}

 Not every $\mathcal{A}$-module $M$ admits a connection in this sense, but it is easy to see that any free $\mathcal{A}$-module does. Of course, arbitrary $\mathcal{A}$-modules do not need to be free. So, in order to obtain a big enough class of modules for which we can guarantee the existence of a Koszul connection, we will make a brief digression on modules of dif\/ferentials and Connes connections (see~\cite{Co85}).

Let $\Omega^1 (\mathcal{A})$ be the $\mathcal{A}$-module def\/ined by the kernel of the multiplication $\mathcal{A} \otimes_{\mathcal{R}} \mathcal{A}\rightarrow \mathcal{A}$. Def\/ine the map $\mathrm{d}:\mathcal{A}\rightarrow \Omega^1 (\mathcal{A})$ by $\mathrm{d}a=1\otimes a -a\otimes 1$, which is a derivation of $\mathcal{A}$ over $\mathcal{R}$ with values into $\Omega^1 (\mathcal{A})$. It is clear from the def\/inition that $\Omega^1 (\mathcal{A})=\mathrm{Span}_{\mathcal{A}}\{ \mathrm{d}f:f\in \mathcal{A}\}$: Since the elements of $\Omega^1 (\mathcal{A})$ lie in the kernel of the multiplication map, if $\sum a_j \otimes b_j \in \Omega^1 (\mathcal{A})$, then $\sum a_jb_j=0$ and therefore
\begin{gather*}
\sum a_j \otimes b_j =\sum (a_j \otimes b_j  -a_jb_j \otimes 1)=\sum a_j\mathrm{d}b_j .
\end{gather*}
 In fact, $\Omega^1 (\mathcal{A})$ is the submodule of $C^1 (\mathrm{Der}(\mathcal{A}),\mathcal{A})$ (the $1$-component of the dif\/ferential algebra $C(\mathrm{Der}(\mathcal{A}),\mathcal{A})$ of Chevalley--Eilenberg cochains of the Lie algebra $\mathrm{Der}(\mathcal{A})$ with values in the $\mathrm{Der}(\mathcal{A})$-module $\mathcal{A}$) generated by the elements $\mathrm{d}f$, $f\in\mathcal{A}$ (see \cite{DM94}). Note, in particular, that this implies $\Omega^1 (\mathcal{A})\subset \mathrm{Der}^* (\mathcal{A})$.
\begin{definition}
Let $M$ be an $\mathcal{A}$-module. A~Connes connection on $M$ is an $\mathcal{A}$-linear map $\delta :M\rightarrow \Omega^1 (\mathcal{A})\otimes_{\mathcal{A}}M$ such that, for all $f\in\mathcal{A}$, $m\in M$,
\begin{gather*}
\delta (fm)=f\delta (m)+\mathrm{d}f \otimes_{\mathcal{A}}m.
\end{gather*}
\end{definition}

\begin{remark}
Connes' def\/inition of a connection (see~\cite{Co85}) actually does not require that $\mathcal{A}$ be a~commutative algebra. The def\/inition goes back to a work by N.~Katz~\cite{Ka70}.
\end{remark}
Starting from a Connes connection, we can obtain a Koszul one. If $X\in \mathrm{Der}(\mathcal{A})$, then we def\/ine a right $\mathcal{A}$-linear pairing $\varphi: \mathrm{Der}(\mathcal{A}) \otimes_{\mathcal{R}}\Omega^1 (\mathcal{A})\rightarrow \mathcal{A}$ by
\begin{gather*}
\Big(X,\sum a_j\mathrm{d}b_j \Big)\mapsto \sum a_j X(b_j ).
\end{gather*}
The Koszul connection $\nabla$ associated to $\delta$ can be constructed as follows: for $X\in \mathrm{Der}(\mathcal{A})$, $\nabla_X \in \mathrm{Hom}_{\mathcal{R}}(M,M)$ is the map given by applying the connection $\delta$ and then contracting the $\Omega^1 (\mathcal{A})$ component with $\varphi$. Thus,
if $m\in M$ is such that $\delta (m)=\sum a_j\mathrm{d}b_j \otimes m_j$, for certain $a_j ,b_j \in \mathcal{A}$ and $m_j \in M$, we have for each $f\in \mathcal{A}$ that $\delta (fm)=f\sum a_j\mathrm{d}b_j \otimes m_j +\mathrm{d}f\otimes m$, and
\begin{gather*}
\nabla_X (fm)  =  \sum \varphi (X,fa_j \mathrm{d}b_j )m_j +\varphi (X,df)m
			   = f\sum a_j X(b_j )m_j +X(f)m \\
\phantom{\nabla_X (fm)}{} = f\nabla_X (m)+X(f)m.
\end{gather*}
A basic result obtained by J.~Cuntz and D.~Quillen (see~\cite{CQ95}) is that Connes connections on an $\mathcal{A}$-module $M$ are in bijective correspondence with $\mathcal{A}$-linear splittings of the natural action $\mathcal{A}\otimes_{\mathcal{R}}M \rightarrow M$. As a consequence, $M$ admits a Connes connection if and only if it is projective.

As said earlier, we will need later on to work with (Koszul) connections, so we need conditions on $\mathcal{A}$ to assure their existence. From what we have seen, these connections exist on any $\mathcal{A}$-module~$M$ which is free or projective. Indeed, note that a free module is always projective, but there are projective modules which are not free. In the literature, there are several well-known conditions on $\mathcal{A}$ guaranteeing the projective character of $M$ (for example, that $\mathcal{A}$ be semi-simple as a ring). When we talk of a connection on $M$, unless otherwise explicitly stated, we will mean that any one of these conditions is satisf\/ied and that the connection is Koszul.

\subsection{Lie algebroids}

\begin{definition}\label{d1}
Let $\mathcal{F}$ be a faithful $\mathcal{A}$-module. A Lie algebroid (on $\mathcal{F}$ over the commutative algebra $\mathcal{A}$) is a triple $(\mathcal{F},[\![ \cdot ,\cdot ]\!],\rho)$, where $[\![ \cdot ,\cdot ]\!]$ is a Lie bracket on $\mathcal{F}$ and $\rho:\mathcal{F}\rightarrow \mathrm{Der}(\mathcal{A})$ is a~morphism of $\mathcal{A}$-modules, called the anchor map, such that
\begin{gather*}
[\![ X ,fY ]\!]=f[\![X,Y]\!]+\rho(X)(f)Y,
\end{gather*}
for all $f\in\mathcal{A}$ and for all $X, Y\in \mathcal{F}$.
\end{definition}

\begin{remark}
Sometimes, the condition that the anchor map be a morphism of Lie algebras is included in the def\/inition of Lie algebroid. However, this fact is a consequence of the conditions in Def\/inition~\ref{d1}, as it has been noted by J.C.~Herz, Y.~Kosmann-Schwarzbach, F.~Magri and J.~Grabowski among others (see~\cite{He53, He53b, KM90, Gra03}).
\end{remark}

\begin{definition}
Let $(\mathcal{F},[\![ \cdot ,\cdot ]\!],\rho)$ and $(\mathcal{F}',[\![ \cdot ,\cdot ]\!]',\rho')$ be Lie algebroids (over the same algebra $\mathcal{A}$ and the same ring $\mathcal{R}$). A morphism of Lie algebroids is a morphism of $\mathcal{A}$-modules $\phi:\mathcal{F}\rightarrow \mathcal{F}'$ such that
\begin{gather*}
\rho'\circ \phi=\rho \qquad \mathrm{and}\qquad \phi([\![ X ,Y ]\!])=[\![ \phi(X) ,\phi(Y) ]\!]',
\end{gather*}
for all $X,Y\in \mathcal{F}$.
\end{definition}
Let us consider some examples. The f\/irst is the classical one.
\begin{example}\label{ec}
Let $M$ be a manifold. Let $E\overset{\pi}{\rightarrow}M$ be a vector bundle over $M$ and
$\mathcal{F}=\Gamma (E)$ the $\mathcal{C}^{\infty}(M)$-module of sections of $E$ (i.e.\ $\mathcal{A}=\mathcal{C}^{\infty}(M)$ and $\mathcal{R}=\mathbb{R}$).
Then, the Lie algebroid structure on $\Gamma(E)$ is def\/ined by a Lie bracket $[\![ \cdot ,\cdot ]\!]$ on $\Gamma (E)$ with an anchor map
\begin{gather*}
q:\Gamma (E)\rightarrow \mathrm{Der}_{\mathbb{R}}\mathcal{C}^{\infty}(M)\cong \Gamma (TM),
\end{gather*}
such that for all $f\in\mathcal{C}^{\infty}(M)$ and for all $X,Y\in \Gamma (E)$:
\begin{enumerate}[1)]\itemsep=0pt
\item $[\![X,fY]\!]=f[\![X,Y]\!]+q(X)(f)Y,$
\item $q(fX+Y)=fq(X)+q(Y).$
\end{enumerate}
In particular, if $E=T^* M$, then the Lie algebroid structure is given by the bracket \eqref{eq1} where the anchor is the Poisson mapping $P$. For $E=TM$ we have the trivial Lie algebroid, where $q =\mathrm{Id}_{TM}$.
\end{example}
\begin{example}
Consider the
$\mathcal{R}$-algebra of dual numbers over $\mathcal{A}$,
\begin{gather*}
\mathcal{A'}=\mathcal{A}[\theta]=\big\{ x+y\theta :x,y\in\mathcal{A},\theta^2=0\big\},
\end{gather*}
with the obvious operations. Clearly, $\mathcal{A'}$ is an $\mathcal{A'}$-module and we can endow it with the Lie algebra structure given by the bracket:
\begin{gather*}
[\![x_1+y_1\theta,x_2+y_2\theta]\!]=(x_1y_2-y_1x_2)\theta
\end{gather*}
for $x_1+y_1\theta,x_2+y_2\theta\in\mathcal{A'}$. Thus $(\mathcal{A'},[\![\cdot,\cdot ]\!], \rho)$ is a Lie algebroid with anchor map
\begin{gather*}
\rho: \ \mathcal{A'} \rightarrow  \mathrm{Der}(\mathcal{A'}), \\
\phantom{\rho: \ }{} \ x+y\theta \longmapsto  \mathrm{ad}_{x}
\end{gather*}
for $x+y\theta\in\mathcal{A'}$. Here $\mathrm{ad}_{x}(x_1+y_1\theta)=[\![x,x_1+y_1\theta]\!]$ is the adjoint map of $x=\mathrm{pr}_1 (x+y\theta )$.
\end{example}
The next example will be relevant in Section~\ref{Transitive}.
\begin{example}
Consider the $\mathcal{A}$-module $\mathrm{Der}(\mathcal{A})\oplus\mathcal{A}$. Denote by $\mathrm{pr}_1$ the projection onto the f\/irst factor, $\mathrm{pr}_1 :\mathrm{Der}(\mathcal{A})\oplus\mathcal{A}\rightarrow \mathrm{Der}(\mathcal{A})$, and def\/ine the following bracket:
\begin{gather*}
[\![(D_1,a_1),(D_2,a_2) ]\!]=([D_1,D_2],D_1(a_2)-D_2(a_1))
\end{gather*}
for $(D_1,a_1),(D_2,a_2)\in \mathrm{Der}(\mathcal{A})\oplus\mathcal{A}$, where $[D_1,D_2]$ is the commutator of endomorphisms. Then, $(\mathrm{Der}(\mathcal{A})\oplus\mathcal{A},[\![\cdot ,\cdot ]\!], \mathrm{pr}_1)$ is a Lie algebroid.
\end{example}

\section{Lie algebroids of Poisson type}\label{LiePoisson}

\begin{definition}
A Poisson algebra $(\mathcal{A},\{\;,\;\})$ is an associative algebra $\mathcal{A}$ together with a Lie bracket which is also a derivation for the product in $\mathcal{A}$, that is, there is an $\mathcal{R}$-bilinear operation
$\{\;,\;\}:\mathcal{A}\times \mathcal{A}\longrightarrow \mathcal{A}$
such that
\begin{enumerate}[1)]\itemsep=0pt
\item $\{f,g\}=-\{g,f\}$ (skew-symmetry),
\item $\{f,\{g,h\}\}+\{g,\{h,f\}\}+\{h,\{f,g\}\}=0$ (Jacobi identity),
\item $\{f,gh\}=\{f,g\}h+g\{f,h\}$ (Leibniz rule),
\end{enumerate}
for all $f,g,h\in\mathcal{A}$.
\end{definition}
If $(\mathcal{A},\{\;,\;\})$ is a Poisson algebra, then we can def\/ine the adjoint map $\mathrm{ad}:\mathcal{A}\longrightarrow \mathrm{Der}(\mathcal{A})$
by
\begin{gather*}
\mathrm{ad}(f)=\{f,\_\},
\end{gather*}
for $f\in\mathcal{A}$. Then, extending the mapping $\mathrm{d}f\longmapsto \mathrm{ad}_f$ by linearity, we get a morphism $\rho: \Omega^1(\mathcal{A})\longrightarrow \mathrm{Der}(\mathcal{A})$ uniquely def\/ined by
\begin{gather}\label{ancla}
\rho(\mathrm{d}f)=\mathrm{ad}_{f},\qquad \forall\,  f\in \mathcal{A}.
\end{gather}
Sometimes (by analogy with Poisson manifolds), $\mathrm{ad}_f$ is referred to as the Hamiltonian vector f\/ield corresponding to $f\in\mathcal{A}$, and denoted by $X_f$ (we will use this notation and terminology later in Section \ref{Poisson}).

Also, given an $\alpha \in\Omega^1(\mathcal{A})$, we can def\/ine $\mathrm{d}\alpha$ through the usual formula
\[
\mathrm{d}\alpha (X,Y)=X(\alpha (Y))-Y(\alpha (X))-\alpha ([X,Y]),
\]
for all $X,Y\in\mathrm{Der}(\mathcal{A})$.
\begin{theorem}\label{rec}
If $(\mathcal{A},\{\;,\;\})$ is a Poisson algebra, then $(\Omega^1(\mathcal{A}),[\;,\;],\rho)$ is a Lie algebroid with anchor map $\rho$ defined by \eqref{ancla} and the Lie bracket
\begin{gather}\label{corch}
[\![\alpha,\beta]\!]=\iota_{\rho(\alpha)}\mathrm{d}\beta-\iota_{\rho(\beta)}\mathrm{d}\alpha+\mathrm{d}(\beta(\rho(\alpha))).
\end{gather}
\end{theorem}

\begin{proof}
All the required properties can be checked directly from the def\/inition of the bracket~\eqref{corch}. To verify the
Jacobi identity, it is easier to consider a system $\{\mathrm{d}f_i\}_{i\in I}$ of generators of $\Omega^1(\mathcal{A})$, and to write the elements $\alpha \in \Omega^1(\mathcal{A})$ as $\alpha=g_i\mathrm{d}f_i$ for some $g_i \in \mathcal{A}$.
\end{proof}

\begin{example}
Let $(M,\{\;,\;\})$ be a Poisson manifold, i.e., $(\mathcal{C}^{\infty}(M),\{\;,\;\})$ is a Poisson algebra. Using Theorem \ref{rec} we have that $\Omega^1(\mathcal{C}^{\infty}(M))=\Omega^1(M)$ is a Lie algebroid with anchor map $\rho =P$ (def\/ined by \eqref{ancla}) and the bracket \eqref{eq1}:
\begin{gather*}
[\![\alpha,\beta]\!]=\iota_{\rho(\alpha)}\mathrm{d}\beta-\iota_{\rho(\beta)}\mathrm{d}\alpha+\mathrm{d}(\beta(\rho(\alpha))),
\end{gather*}
for all $\alpha,\beta \in \Omega^1(M)$. The proof is based on the fact that there exists a f\/inite subset $\{ g_1 ,\dots ,g_k \}\!\subset {C}^{\infty}(M)$ (with $k\!\leq\! 2\dim M\!+\!1$) such that the ${C}^{\infty}(M)$-module $\Omega^1(M)$ is spanned by
$\{ \mathrm{d}g_1 ,{\dots},\mathrm{d}g_k \}$. This, in turn, is a consequence of Whitney's embedding theorem and the fact that
the sheaf of germs of smooth functions is soft (see~\cite{AG90}).
\end{example}

Theorem \ref{rec} tells us that given a Poisson algebra we have a Lie algebroid canonically associated to it. We are interested now in the reciprocal: When does a Lie algebroid $(\Omega^1(\mathcal{A}),[\![\cdot ,\cdot ]\!],\rho)$ determine a Poisson structure on $\mathcal{A}$?.

Given a Lie algebroid $(\Omega^1(\mathcal{A}),[\![\cdot ,\cdot ]\!],\rho)$, for any $f,g\in\mathcal{A}$ we can def\/ine
\begin{gather}\label{cp}
\{f,g\}=(\mathrm{d}g)(\rho(\mathrm{d}f))=(\rho(\mathrm{d}f))(g).
\end{gather}
The bracket $\{\;,\;\}$ def\/ined in this way is clearly $\mathcal{R}$-bilinear and satisf\/ies the Leibniz rule
\begin{gather*}
\{f,gh\} = (\rho(df))(gh)=g\rho(df)(h)+\rho(df)(g)h=g\{f,h\}+\{f,g\}h,
\end{gather*}
 for all $f,g,h\in \mathcal{A}$.

Thus, in order to get a Poisson structure on $\mathcal{A}$ we only need to take care of the skew-symmetry and the Jacobi identity.
\begin{definition}
The anchor map $\rho$ of a Lie algebroid $(\Omega^1(\mathcal{A}),[\![\cdot ,\cdot ]\!],\rho)$ is said to be skew-symmet\-ric~if
\begin{gather*}
\alpha(\rho(\beta))=-\beta(\rho(\alpha))
\end{gather*}
for all $\alpha, \beta\in \Omega^1(\mathcal{A})$.
\end{definition}
If we assume that the anchor map is skew-symmetric, then the new operation def\/ined by \eqref{cp} is also skew-symmetric. Now,  let us turn our attention to the Jacobi identity.
\begin{theorem}\label{tsn}
Let $(\Omega^1(\mathcal{A}),[\![\cdot ,\cdot ]\!],\rho)$ be a Lie algebroid, and define $Q\in \Lambda^2(\mathrm{Der}(\mathcal{A}))$ by $Q(\mathrm{d}f,\mathrm{d}g)=\mathrm{d}g(\rho(\mathrm{d}f))$ for $f,g\in\mathcal{A}$. The following  conditions are equivalent:
\begin{enumerate}[$(i)$]\itemsep=0pt
\item the bracket \label{ij} $\{\;,\;\}$ defined by $\{f,g\}=\rho(\mathrm{d}f)(g)$ satisfies the Jacobi identity,
\item \label{sn} $[Q,Q]_{\rm SN}=0$,
\item \label{qm}$[\rho(\mathrm{d}f),\rho(\mathrm{d}g)]=\rho(\mathrm{d}\{f,g\})$.
\end{enumerate}
Here, $[\cdot ,\cdot ]_{\rm SN}$ denotes the Schouten--Nijenhuis bracket on multiderivations \emph{\cite{KSM93}}.
\end{theorem}

\begin{proof}
It is based on the independent computation of both sides of the equation $\frac{1}{2}[\mathcal{L}_Q,\iota_Q]=\frac{1}{2}\iota_{[Q,Q]_{\rm SN}}$, where on the left-hand side we have the graded commutator of operators on the graded algebra $\Lambda(\mathcal{A})$. Note that due to the $\mathcal{A}$-linearity, we only need to apply these operators to basis elements.
\end{proof}

This result motivates the following def\/inition.
\begin{definition}
A Lie algebroid $(\Omega^1(\mathcal{A}),[\![\cdot ,\cdot ]\!],\rho)$ is of Poisson type if:
\begin{enumerate}[1)]\itemsep=0pt
\item the anchor $\rho$ is skew-symmetric,
\item one of the equivalent conditions \eqref{ij}--\eqref{qm} in Theorem~\ref{tsn} holds.
\end{enumerate}
\end{definition}

 In other words, a Lie algebroid $(\Omega^1(\mathcal{A}),[\![\cdot ,\cdot ]\!] ,\rho)$ is of Poisson type if it determines a Poisson structure on the algebra $\mathcal{A}$.

Another important issue for us is to determine the form of the bracket of a Lie algebroid of Poisson type. We know that the classical example of Poisson manifolds leads to brackets of the type \eqref{eq1}. The following result characterizes the class of such algebroids.

\begin{proposition}\label{c}\sloppy
If $(\Omega^1(\mathcal{A}),[\![\cdot ,\cdot ]\!],\rho)$ is a Lie algebroid for which $\rho$ is skew-symmetric and $[\![\mathrm{d}f,\mathrm{d}g]\!]=\mathrm{d}\{f,g\}$, then the bracket of the Lie algebroid is of the form
\begin{gather*}
[\![\alpha,\beta]\!]=\iota_{\rho(\alpha)}\mathrm{d}\beta-\iota_{\rho(\beta)}\mathrm{d}\alpha+\mathrm{d}(\beta(\rho(\alpha))).
\end{gather*}
\end{proposition}

\begin{proof}
Note f\/irst that, if $\alpha=\mathrm{d}f$ and $\beta=\mathrm{d}g$ for some $f,g,\in\mathcal{A}$, then
\begin{gather*}
\iota_{\rho(\alpha)}d\beta -\iota_{\rho(\beta)}\mathrm{d}\alpha+\mathrm{d}(\beta(\rho(\alpha)))=\mathrm{d}(\mathrm{d}g(\rho(\mathrm{d}f)))=\mathrm{d}\{f,g\}=[\![\mathrm{d}f,\mathrm{d}g]\!]=[\![\alpha,\beta]\!].
\end{gather*}
Since every element of $\Omega^1(\mathcal{A})$ is a linear combination of elements of the form $\mathrm{d}f_i$ ($f_i\in\mathcal{A}$), it is enough to prove the statement for $\alpha=f_1\mathrm{d}g_1$ and $\beta=f_2\mathrm{d}g_2$,
which can be done by a direct computation.
\end{proof}
The hypothesis of the Proposition \ref{c} can be reformulated in an alternative way.
\begin{proposition}\label{eq}
If $(\Omega^1(\mathcal{A}),[\![\cdot ,\cdot ]\!],\rho)$ is a Lie algebroid with anchor $\rho$ skew-symmetric, then the following assertions are equivalent:
\begin{enumerate}[$(i)$]\itemsep=0pt
\item \label{1} $[\![\mathrm{d}f,\mathrm{d}g ]\!]=\mathrm{d}(\mathrm{d}g(\rho(\mathrm{d}f)))=\mathrm{d}\{f,g\}$ for all $f,g\in\mathcal{A}$,
\item \label{2} $\mathrm{d}\alpha=0=\mathrm{d}\beta$ implies $\mathrm{d}[\![\alpha,\beta ]\!]=0$ for $\alpha, \beta \in \Omega^1(\mathcal{A})$.
\end{enumerate}
\end{proposition}
\begin{proof}
First, let us assume that item \eqref{1} holds. By Proposition \ref{c} we know that the bracket is of the form
\begin{gather*}
[\![\alpha,\beta]\!]=\iota_{\rho(\alpha)}\mathrm{d}\beta-\iota_{\rho(\beta)}\mathrm{d}\alpha+\mathrm{d}(\beta(\rho(\alpha))),
\end{gather*}
and hence condition \eqref{2} holds,
\begin{gather*}
\mathrm{d}[\![\alpha,\beta]\!]=\mathrm{d}(\iota_{\rho(\alpha)}\mathrm{d}\beta-\iota_{\rho(\beta)}\mathrm{d}\alpha+\mathrm{d}(\beta(\rho(\alpha))))=\mathrm{d}^2(\beta(\rho(\alpha)))=0,
\end{gather*}
whenever $\mathrm{d}\alpha=0=\mathrm{d}\beta$.

Now let us assume that condition \eqref{2} holds. Def\/ine $C(\alpha,\beta)=[\![\alpha,\beta]\!]-\mathrm{d}(\beta(\rho(\alpha)))$ for $\alpha,\beta\in \Omega^1(\mathcal{A})$.
Notice that $C$ is skew-symmetric and
\begin{gather*}
\mathrm{d}C(\alpha,\beta)=\mathrm{d}([\![\alpha,\beta]\!])-\mathrm{d}^2(\beta(\rho(\alpha)))=0,
\end{gather*}
for any closed $\alpha$ and $\beta$.

Let us evaluate the following expression
\begin{gather*}
C(\mathrm{d}f,h\mathrm{d}h)=[\![\mathrm{d}f,h\mathrm{d}h]\!]-\mathrm{d}(h\mathrm{d}h(\rho(\mathrm{d}f)))\\
\phantom{C(\mathrm{d}f,h\mathrm{d}h)}{}
=h[\![\mathrm{d}f,\mathrm{d}h]\!]+\rho(\mathrm{d}f)(\mathrm{d}h)\mathrm{d}h-\mathrm{d}h(\mathrm{d}h(\rho(\mathrm{d}f)))-h\mathrm{d}(h\mathrm{d}h(\rho(\mathrm{d}f)))\\
\phantom{C(\mathrm{d}f,h\mathrm{d}h)}{}
=h[\![\mathrm{d}f,\mathrm{d}h]\!]-h\mathrm{d}(h\mathrm{d}h(\rho(\mathrm{d}f)))
=hC(\mathrm{d}f,\mathrm{d}h).
\end{gather*}
Now, applying the operator $\mathrm{d}$ and taking into account that $h\mathrm{d}h$ is closed, we get
\begin{gather*}
0=\mathrm{d}C(\mathrm{d}f,h\mathrm{d}h)=\mathrm{d}(hC(\mathrm{d}f,\mathrm{d}h))=\mathrm{d}h\wedge C(\mathrm{d}f,\mathrm{d}h)+h\mathrm{d}C(\mathrm{d}f,\mathrm{d}h)=\mathrm{d}h\wedge C(\mathrm{d}f,\mathrm{d}h).
\end{gather*}
Interchanging the roles of $f$ and $h$, we have
\begin{gather*}
0=\mathrm{d}C(\mathrm{d}h,f\mathrm{d}f)=\mathrm{d}f\wedge C(\mathrm{d}h,\mathrm{d}f)=-\mathrm{d}f\wedge C(\mathrm{d}f,\mathrm{d}h).
\end{gather*}
These relations imply that $\mathrm{d}f$, $\mathrm{d}h$ and $C(\mathrm{d}f,\mathrm{d}h)$ are linearly dependent for arbitrary $f$ and $h$. In particular, if $\mathrm{d}f$ and $\mathrm{d}h$ are linearly independent, then $C(\mathrm{d}f,\mathrm{d}h)=0$, and hence
\begin{gather*}
[\![\mathrm{d}f,\mathrm{d}h]\!]=\mathrm{d}(\mathrm{d}h(\rho(\mathrm{d}f)))=\mathrm{d}\{f,h\}.\tag*{\qed}
\end{gather*}
\renewcommand{\qed}{}
\end{proof}

\begin{remark}
It has been shown that a Poisson bracket $\{\cdot,\cdot\}$ on $M$ def\/ines a graded Poisson bracket $[\![\cdot,\cdot]\!]$ in $\Omega(M)$ (see~\cite{Kz85}). Its values on smooth functions and on exact 1-forms are given by $[\![f,g]\!]=0$, $[\![\mathrm{d}f,g]\!]=\{f,g\}$ and $[\![\mathrm{d}f,\mathrm{d}g]\!]=\mathrm{d}\{f,g\}$, respectively. Then, the bracket can be extended to all dif\/ferential forms by the graded Leibniz rule. For this Poisson bracket the exterior dif\/ferential is a Poisson derivation, that is,
\[
\mathrm{d}[\![\alpha,\beta]\!]=[\![\mathrm{d}\alpha,\beta]\!]+(-1)^{|\alpha|-1}[\![\alpha,\mathrm{d}\beta]\!]
\]
for any $\alpha,\beta\in\Omega(M)$. Such graded Poisson brackets were called dif\/ferential in~\cite{Be95}.
\end{remark}

Let us summarize these results.
\begin{theorem}\label{i}
Let $(\Omega^1(\mathcal{A}),[\![\cdot ,\cdot ]\!],\rho)$ be a Lie algebroid. Then, there is a Poisson algebra structure $\{\;,\;\}:\mathcal{A}\times \mathcal{A}\rightarrow \mathcal{A}$ on $\mathcal{A}$ such that
\begin{gather*}
\rho(\mathrm{d}f)=\mathrm{ad}_f ,\qquad \forall\, f\in\mathcal{A},
\end{gather*}
if and only if:
\begin{enumerate}[$(a)$]\itemsep=0pt
\item $\rho$ is skew-symmetric \label{ca}
\item One of the following conditions holds:
	\begin{enumerate}[$1)$]\itemsep=0pt
	\item \label{cb} $\mathrm{d}\alpha=0=\mathrm{d}\beta$ implies $\mathrm{d}[\![\alpha,\beta]\!]=0$;
	\item $[\![\mathrm{d}f,\mathrm{d}g]\!]=\mathrm{d}(\mathrm{d}g(\rho(\mathrm{d}f)))$ for all $f,g\in\mathcal{A}$.
	\end{enumerate}
\end{enumerate}
Under these conditions, the Lie bracket $[\![\cdot ,\cdot ]\!]$ is reconstructed from $\rho$ by the formula
\begin{gather*}
[\![\alpha ,\beta ]\!]=\iota_{\rho(\alpha)}\mathrm{d}\beta-\iota_{\rho(\beta)}\mathrm{d}\alpha+\mathrm{d}(\beta(\rho(\alpha))).
\end{gather*}
\end{theorem}

Of course, the basic example of this situation is the cotangent Lie algebroid of a symplectic manifold. A Lie algebroid structure $(\Omega^1(M),[\![\cdot ,\cdot ]\!],\rho )$ induces a Poisson bracket on $\mathcal{C}^{\infty}(M)$ if and only if $\rho$ is skew-symmetric and, whenever $\mathrm{d}\alpha=0=\mathrm{d}\beta$, then $\mathrm{d}[\![\alpha,\beta]\!]=0$.
For this kind of examples, the property $[\![\mathrm{d}f,\mathrm{d}g]\!]=\mathrm{d}\{f,g\}$ follows directly from the injectivity of the anchor map and the fact that it is a Lie algebra morphism: $\rho [\![\mathrm{d}f,\mathrm{d}g]\!]=[\![\rho (\mathrm{d}f),\rho (\mathrm{d}g)]\!]=\rho (\mathrm{d}\{f,g\})$. Let us consider the following example, where the anchor map is not injective but the property $[\![\mathrm{d}f,\mathrm{d}g]\!]=\mathrm{d}\{f,g\}$ still holds.
\begin{example}
Let $\mathcal{R}=\mathbb{R}$ and $\mathcal{A}=\mathbb{R}[x^1,x^2,x^3]$. Then,
$\mathrm{Der}(\mathcal{A})=\mathrm{Span}\left\{\partial_1,\partial_2,\partial_2 \right\}$ and $\Omega^{1}(\mathcal{A})=\mathrm{Span}\{\mathrm{d}x^1,\mathrm{d}x^2,\mathrm{d}x^3\}$. Def\/ine the following bracket
\begin{gather*}
[\![p_i\mathrm{d}x^i,q_j\mathrm{d}x^j]\!]=\left(\left(
p_1\left(\partial_2+\partial_3\right)+
p_2\left(-\partial_1+\partial_3\right)+
p_3\left(-\partial_1-\partial_2\right)\right)(q_i)\right.\\
\left. \hphantom{[\![p_i\mathrm{d}x^i,q_j\mathrm{d}x^j]\!]=}{}
-\left(q_1\left(\partial_2+\partial_3\right)+
q_2\left(-\partial_1+\partial_3\right)+
q_3\left(-\partial_1-\partial_2\right)\right)(p_i)\right) \mathrm{d}x^i ,
\end{gather*}
and the anchor map as
\begin{gather*}
\rho: \ \Omega^1(\mathcal{A}) \rightarrow \mathrm{Der}(\mathcal{A}), \\
\phantom{\rho: \ }{} \ p_i\mathrm{d}x^i \mapsto -(p_2+p_3)\partial_1+(p_1-p_3)\partial_2+(p_1+p_2)\partial_3.
 \end{gather*}
Note that the matrix representation of $\rho$ relative to the given basis in $\mathrm{Der}(\mathcal{A})$ and $\Omega^{1}(\mathcal{A})$ is
\begin{gather*}
\rho =
\begin{pmatrix}
0 & 1 & 1 \\
-1 & 0 & 1 \\
-1 & -1 & 0
\end{pmatrix}.
\end{gather*}
Therefore, $\rho$ is skew-symmetric and of rank $2$ ($\rho$ is not injective).
A long but straightforward computation shows that $(\Omega^1(\mathcal{A}),]\![\cdot ,\cdot ]\!],\rho)$ is a Lie algebroid.

Let us show that this Lie algebroid is of Poisson type by checking the property $[\![df,dg]\!]=d\{f,g\}$. We have for $p,q\in\mathcal{A}$:
\begin{gather*}
\{p,q\} = \mathrm{d}(\rho (\mathrm{d}p)(q))= \mathrm{d}\left( \rho \left( \partial_i p\mathrm{d}x^i \right) (q)\right) \\
\phantom{\{p,q\}}{} = \mathrm{d}\left( -(\partial_2p+\partial_3p)\partial_1q
	   +(\partial_1p-\partial_3p)\partial_2q
	   +(\partial_1p+\partial_2p)\partial_3q \right).
\end{gather*}
The $\mathrm{d}x^1$ factor in the expansion of this expression (the other cases are similar) is
\begin{gather*}
 - (\partial_2p+\partial_3p)\partial^2_{11}q +(\partial_1p-\partial_3p)\partial^2_{12}q +(\partial_1p+\partial_2p)\partial^2_{13}q \\
\qquad{}- (\partial^2_{12}p+\partial^2_{13}p)\partial_1q+(\partial^2_{11}p-\partial^2_{13}p)\partial_2q + (\partial^2_{11}p+\partial^2_{12}p)\partial_3q.
\end{gather*}
On the other hand, the Lie algebroid bracket is
\begin{gather*}
[\![\mathrm{d}p,\mathrm{d}q]\!] =  [\![\partial_i p\mathrm{d}x^i,\partial_j q\mathrm{d}x^j]\!]
							    =  \rho(\mathrm{d}p)(\partial_k q)-\rho(\mathrm{d}q)(\partial_k p))\mathrm{d}x^k.
\end{gather*}
For $k=1$, we compute the coef\/f\/icient of $\mathrm{d}x^1$:
\begin{gather*}
 - (\partial_2p+\partial_3p)\partial^2_{11}q +(\partial_1p-\partial_3p)\partial^2_{12}q +(\partial_1p+\partial_2p)\partial^2_{13}q \\
\qquad{} + (\partial_2q+\partial_3q)\partial^2_{11}p -(\partial_1q-\partial_3q)\partial^2_{12}p -(\partial_1q+\partial_2q)\partial^2_{13}p,
\end{gather*}
which is the same as above. Thus, $(\Omega^1(\mathcal{A}),[\![\cdot ,\cdot ]\!],\rho)$ is a Lie algebroid of Poisson type and its bracket is just given by formula~\eqref{corch}.
\end{example}
\begin{remark}
The anchor map of the cotangent Lie algebroid is not injective in general. It is only true in the symplectic case.
\end{remark}

Finally, let us consider an example of Lie algebroid of Poisson type whose bracket does not have the form~\eqref{corch}.
\begin{example}
Let $\mathcal{R}=\mathbb{R}$ and $\mathcal{A}=\mathbb{R}[x^1,x^2,x^3]$. Def\/ine the structure of an Abelian Lie algebra with generators $\{\mathrm{d}x^1,\mathrm{d}x^2,\mathrm{d}x^3\}$,
\begin{gather*}
[\![\mathrm{d}x^i,\mathrm{d}x^j]\!]=0,\qquad i,j\in\{1,2,3\}.
\end{gather*}
Extending by $\mathcal{R}$-bilinearity and the Leibniz identity gives
\begin{gather*}
[\![p_i\mathrm{d}x^i,q_j\mathrm{d}x^j]\!]=p_i\rho(\mathrm{d}x^i)(q_j)\mathrm{d}x^j-q_j\rho(\mathrm{d}x^j)(p_i)\mathrm{d}x^i.
\end{gather*}
Next, let us think of the following skew-symmetric morphism of  $\mathcal{A}$-modules as the anchor map:
\begin{gather*}
\rho: \ \Omega^1(\mathcal{A}) \rightarrow \mathrm{Der}(\mathcal{A}),\\
\phantom{\rho: \ }{} \ \mathrm{d}x^1 \mapsto x^3\partial_2,\\
\phantom{\rho: \ }{} \ \mathrm{d}x^2 \mapsto -x^3\partial_1,\\
\phantom{\rho: \ }{} \ \mathrm{d}x^3 \mapsto 0.
\end{gather*}
With these def\/initions the bracket in $\Omega^1(\mathcal{A})$ can also be expressed in a form suitable for explicit computations, as
 \begin{gather*}
 [\![p_i\mathrm{d}x^i,p_j\mathrm{d}x^j]\!]=x^3((p\times \nabla)_3q_i-(q\times \nabla)_3p_i)\mathrm{d}x^i,
 \end{gather*}
where an element $p\in\mathcal{A}$ is viewed as a vector $p=(p_1,p_2,p_3)$, $\nabla =(\partial_1,\partial_2,\partial_3)$ and the subindex~$3$ denotes the third component of the cross product $p\times \nabla$.

It can be checked that $(\Omega^1(\mathcal{A}),[\![\cdot,\cdot]\!],\rho)$ is a Lie algebroid of Poisson type. However, if $p=2x^1+x^2$ and $q=x^1+x^2$, then
 \begin{gather*}
 [\![\mathrm{d}p,\mathrm{d}q]\!]=[\![2\mathrm{d}x^1+\mathrm{d}x^2,\mathrm{d}x^1+\mathrm{d}x^2]\!]=0,
 \end{gather*}
whereas
 \begin{gather*}
 \mathrm{d}\{p,q\}=\mathrm{d}(\mathrm{d}q(\rho(\mathrm{d}p))) = \mathrm{d}((\mathrm{d}x^1+\mathrm{d}x^2)\rho(2\mathrm{d}x^1+\mathrm{d}x^2))\\
\phantom{\mathrm{d}\{p,q\}=\mathrm{d}(\mathrm{d}q(\rho(\mathrm{d}p)))}{}
= \mathrm{d}((\mathrm{d}x^1+\mathrm{d}x^2)(2x^3\partial_2 -x^3\partial_1))=\mathrm{d}x^3.
 \end{gather*}
 \end{example}

\section{Transitive Lie algebroids}\label{Transitive}

To motivate the def\/inition of transitive Lie algebroids let us consider for a moment the geometric example of a Lie algebroid $(E,[\![\cdot ,\cdot ]\!],q)$, where $E\rightarrow M$ is a vector bundle over a manifold $M$ (recall Example \ref{ec}). If the anchor map $q:\Gamma(E)\rightarrow \Gamma(TM)$ is an epimorphism, the algebroid $(E,[\![\cdot ,\cdot ]\!],q)$ is said to be transitive. In this case, it is possible to construct the so-called Atiyah sequence of the algebroid, which is the short exact sequence
\begin{gather*}
\begin{diagram}
\node{\mathfrak{g}} \arrow{e,t,J}{j} \node{E} \arrow{e,t,A}{q} \node{TM}
\end{diagram}
\end{gather*}
where $\mathfrak{g}=\mathrm{Ker}\,q$. Thus, the existence of a section for the anchor $q$ (equivalently, a linear connection on $E$) implies that, locally, $E=TM\oplus \mathfrak{g}$. Note also that the f\/ibre of the bundle $\mathfrak{g}$ over the point $x\in M$, $\mathfrak{g}_x$ is a Lie algebra (called the isotropy Lie algebra of the algebroid $E$ at $x\in M$) with the bracket given, for $\alpha,\beta\in\mathfrak{g}_x$, by
\begin{gather*}
[\alpha ,\beta ]=[\![X,Y]\!],
\end{gather*}
where $X,Y\in\Gamma(E)$ are any sections such that $X(x)=\alpha$ and $Y(x)=\beta$.
\begin{definition}
A Lie algebroid $(\mathcal{F},[\![\cdot ,\cdot ]\!],\rho)$ is transitive if there exist a short exact sequence of $\mathcal{A}$-modules:
\begin{gather*}
\begin{diagram}
\node{V} \arrow{e,t,J}{j} \node{\mathcal{F}} \arrow{e,t,A}{\rho} \node{\mathrm{Der}(\mathcal{A})}
\end{diagram}
\end{gather*}
\end{definition}
\begin{remark}
There is a more general notion, the extension of a Lie--Rinehart algebra, that generalizes the transitivity condition for a geometric Lie algebroid (see \cite{Hue99}).
\end{remark}

\subsection{Transitive algebroids induced by connections}
Let $V$ be a unitary $\mathcal{A}$-module and $\nabla$ a connection on $V$.
\begin{definition}
The curvature of $\nabla$ is the mapping $C_{\nabla}:\mathrm{Der}(\mathcal{A})\times\mathrm{Der}(\mathcal{A})\rightarrow \mathrm{Hom}(V,V)$ def\/ined by
\begin{gather*}
C_{\nabla}(X,Y)=[\nabla_{X},\nabla_{Y}]-\nabla_{[X,Y]}.
\end{gather*}
\end{definition}

\begin{definition}
If $V$ is endowed with a Lie algebra structure $[\cdot,\cdot]_V$, then a connection $\nabla$ is said to be a Lie connection if
\begin{gather*}
\nabla_X[v_1,v_2]_V=[\nabla_Xv_1,v_2]_V+[v_1,\nabla_Xv_2]_V,
\end{gather*}
for all $v_1,v_2\in V$ and for all $X\in \mathrm{Der}(\mathcal{A})$.
\end{definition}

We have the following result, which gives us a description of transitive Lie algebroids in terms of Lie connections (see, also \cite{Mc95, Mc05}).

\begin{theorem}\label{tat}
Let $V$ be an $\mathcal{A}$-module endowed with a Lie algebra structure $[\;,\;]_V$ which is $\mathcal{A}$-linear $($i.e., $[fv_1,v_2]_V=f[v_1,v_2]_V$ $\forall\, f\in\mathcal{A}$, $v_1,v_2\in V)$. Let $\nabla$ be a Lie connection on $V$.
If there exists a $2$-form $\mathcal{B}\in \Omega^2 (\mathcal{A};V)$ with values in $V$, such that, for any $X_1, X_2, X_3 \in \mathrm{Der}(\mathcal{A})$ and $v\in V$ the following conditions hold:
\begin{enumerate}[$(a)$]\itemsep=0pt
\item $[\mathcal{B}(X_1,X_2),v]_V=C_{\nabla}(X_1,X_2)(v)$, \label{tat1}
\item the cyclic sum $\circlearrowleft\left( \nabla_{X_1}(\mathcal{B}(X_2,X_3))-\mathcal{B}([X_1,X_2],X_3)\right) =0$, \label{tat2}
\end{enumerate}
then $(\mathrm{Der}(\mathcal{A})\oplus V,[\![\cdot ,\cdot ]\!],\rho)$ is a transitive Lie algebroid
with anchor map $\rho=\mathrm{pr}_1$ and bracket
\begin{gather}\label{cl}
[\![(X_1,v_1),(X_2,v_2)]\!]=([X_1,X_2],[v_1,v_2]_V+\nabla_{X_1}v_2-\nabla_{X_2}v_1-\mathcal{B}(X_1,X_2)).
\end{gather}
Moreover,
\begin{gather*}
\iota_2(\mathcal{B}(X,Y))=\iota_1([X,Y])-[\![\iota_1(X),\iota_1(Y)]\!],
\end{gather*}
for $X,Y\in \mathrm{Der}(\mathcal{A})$. Here $\iota_1 :\mathrm{Der}(\mathcal{A})\rightarrow \mathrm{Der}(\mathcal{A})\oplus V$ and $\iota_2 :V\rightarrow \mathrm{Der}(\mathcal{A})\oplus V$ are the inclusion maps.
\end{theorem}

\begin{proof}
It is clear that $ \mathrm{Der}(\mathcal{A})\oplus V$ is an $\mathcal{A}$-module. We must show that the bracket def\/ined by~\eqref{cl} is Lie. The skew-symmetry and the $\mathcal{R}$-bilinearity are immediate, while the Jacobi identity and the Leibniz rule are checked by lengthy but straightforward computations. Since $\rho=\mathrm{pr}_1$ is clearly $\mathcal{A}$-linear, we have that $\mathrm{Der}(\mathcal{A})\oplus V$ is a Lie algebroid. The transitivity is obvious in view of the sequence
\begin{gather*}
\begin{diagram}
\node{V} \arrow{e,t,J}{\iota_2} \node{\mathrm{Der}(\mathcal{A})\oplus V} \arrow{e,t,A}{\iota_1} \node{\mathrm{Der}(\mathcal{A}).}
\end{diagram}
\end{gather*}
Moreover, a direct computation shows that
\begin{gather*}
\iota_1[X,Y]-[\![\iota_1(X),\iota_1(Y)]\!] = ([X,Y],0)-[\![(X,0),(Y,0)]\!]
 = ([X,Y],0)-([X,Y],-\mathcal{B}(X,Y)) \\
\phantom{\iota_1[X,Y]-[\![\iota_1(X),\iota_1(Y)]\!]}{}
 =  (0,\mathcal{B}(X,Y))=\iota_2(\mathcal{B}(X,Y)).\tag*{\qed}
 \end{gather*}
\renewcommand{\qed}{}
\end{proof}

\begin{remark}
The conditions \eqref{tat1} and~\eqref{tat2} in this theorem have the following interpretation. Condition \eqref{tat1} states that the curvature of the connection $\nabla$ is given by the composition $C_{\nabla}=\mathrm{ad}^R \circ \mathcal{B}$, where $\mathrm{ad}^R :V\rightarrow V$ is the right adjoint with respect to the bracket on~$V$. On the other hand, \eqref{tat2} expresses the Bianchi identity for a connection with torsion (see also \cite{Mc95,Mc05}).
\end{remark}

\subsection[Parametrization of transitive algebroids on $\mathrm{Der}\mathcal{A}\oplus V$]{Parametrization of transitive algebroids on $\boldsymbol{\mathrm{Der}\mathcal{A}\oplus V}$}

The following result states that the converse of Theorem \ref{tat} is also true.

\begin{theorem}\label{tatr}
Let $(\mathrm{Der}(\mathcal{A})\oplus V,[\![\cdot ,\cdot ]\!],\rho)$ be a transitive Lie algebroid with
$\rho=\mathrm{pr}_1$. Then:
\begin{enumerate}[$(i)$]\itemsep=0pt
\item The bracket on $V$, defined for $v_1,v_2\in V$ by
\begin{gather*}
[v_1,v_2]_V=\mathrm{pr}_2([\![\iota_2(v_1),\iota_2(v_2)]\!]),
\end{gather*}
is an $\mathcal{A}$-linear Lie bracket.
\item The mapping $\nabla :\mathrm{Der}(\mathcal{A})\rightarrow \mathrm{Hom}(V,V)$ given by
\begin{gather*}
\nabla_X(v)=\mathrm{pr}_2([\![\iota_1(X),\iota_2(v)]\!]),
\end{gather*}
for $v\in V$ and $X\in \mathrm{Der}(\mathcal{A})$, is a Lie connection on $V$.
\item The mapping $\mathcal{B}:\mathrm{Der}(\mathcal{A})\times \mathrm{Der}(\mathcal{A})\rightarrow V$ defined by
\begin{gather*}
\mathcal{B}(X,Y)=\mathrm{pr}_2([\![\iota_1(X),\iota_1(Y)]\!]-\iota_1([X,Y])),
\end{gather*}
is $\mathcal{A}$-bilinear, skew-symmetric and it satisfies conditions \eqref{tat1} and \eqref{tat2} of Theorem {\rm \ref{tat}}.
\end{enumerate}
\end{theorem}

\begin{proof}
The proof is a lengthy but straightforward computation developing the def\/initions.
\end{proof}

As a consequence of Theorems \ref{tat} and \ref{tatr}, we have the following.
\begin{corollary}
Let $V$ be an $\mathcal{A}$-module. There is a one to one correspondence between transitive Lie algebroids $(\mathrm{Der}(\mathcal{A})\oplus V,[\![\;,\;]\!],\rho)$ with anchor $\rho=\mathrm{pr}_1$ and triples $([\cdot ,\cdot ]_V ,\nabla ,\mathcal{B})$ consisting of:
\begin{enumerate}[$(i)$]\itemsep=0pt
\item an $\mathcal{A}$-linear Lie bracket $[\cdot ,\cdot ]_V$ in $V$;
\item a Lie connection $\nabla$ on $V$;
\item a $2$-form $\mathcal{B}\in\Omega^2 (\mathcal{A};V)$
satisfying conditions \eqref{tat1} and \eqref{tat2} of Theorem {\rm \ref{tat}}.
\end{enumerate}
\end{corollary}

\subsection{Algebroid connections}

In this subsection, we f\/irst generalize Theorem \ref{tatr} to the case of a transitive Lie algebroid which is not necessarily of the form $\mathrm{Der}(\mathcal{A})\oplus V$, but we require that it be endowed with a Lie algebroid connection (see \cite{Ku91,Mc05}), considered as a section of the anchor map. Then, in Theorem \ref{iso} we construct a Lie algebroid structure on $\mathrm{Der}(\mathcal{A})\oplus V$ which is isomorphic to the given algebroid.
\begin{definition}
Let $(\mathcal{F},[\![\cdot ,\cdot ]\!],\rho)$ be a Lie algebroid. A morphism of $\mathcal{A}$-modules $\gamma:\mathrm{Der}(\mathcal{A})\rightarrow \mathcal{F}$ is called a Lie algebroid connection (on $\mathcal{F}$) if it is a section of $\rho$, that is, $\rho\circ\gamma=\mathrm{Id}_{\mathrm{Der}(\mathcal{A})}$.
\end{definition}

As stated in Section~\ref{subsec21}, we will assume some condition guaranteeing the existence of a~connection on $(\mathcal{F},[\![\cdot ,\cdot ]\!],\rho)$, to be precise, we will assume that $\mathcal{F}$ is free. It is known that the existence of a section $\gamma$ of $\rho$ in the short exact sequence
\begin{gather*}
\begin{diagram}
\node{V} \arrow{e,t,J}{\iota} \node{\mathcal{F}} \arrow{e,t,A}{\rho} \node{\mathrm{Der}(\mathcal{A})}
\end{diagram}
\end{gather*}
implies that $\mathcal{F}\cong V\oplus \mathrm{Der}(\mathcal{A})$.

\begin{proposition}\label{g}
Let
\begin{gather}\label{exact}
\begin{diagram}
\node{V} \arrow{e,t,J}{\iota} \node{\mathcal{F}} \arrow{e,t,A}{\rho} \node{\mathrm{Der}(\mathcal{A})}
\end{diagram}
\end{gather}
be a transitive Lie algebroid and $\gamma:\mathrm{Der}(\mathcal{A})\rightarrow \mathcal{F}$ a Lie algebroid connection for $\mathcal{F}$.
Then:
\begin{enumerate}[$(i)$]\itemsep=0pt
\item The bracket $[\cdot ,\cdot ]_V$ defined for $v_1,v_2\in V$ by
\begin{gather*}
[v_1,v_2]_V=[\![\iota(v_1),\iota(v_2)]\!],
\end{gather*}
is an $\mathcal{A}$-linear Lie bracket.
\item The mapping $\nabla :\mathrm{Der}(\mathcal{A})\rightarrow \mathrm{Hom}(V,V)$ given by
\begin{gather*}
\nabla_X(v)=[\![\gamma(X),\iota(v)]\!],
\end{gather*}
is a Lie connection on $V$ $($called the adjoint connection$)$.
\item The mapping  $\mathcal{B}:\mathrm{Der}(\mathcal{A})\times \mathrm{Der}(\mathcal{A})\rightarrow V$, defined by
\begin{gather*}
\mathcal{B}(X,Y)=[\![\gamma(X),\gamma(Y)]\!]-\gamma([X,Y]),
\end{gather*}
is $\mathcal{A}$-bilinear, skew-symmetric and satisfies the conditions \eqref{tat1} and \eqref{tat2} of Theorem {\rm \ref{tat}}.
\end{enumerate}
\end{proposition}

\begin{proof}
First, note that the exactness of \eqref{exact} allows us to identify $V$ with its image under $\iota$ in~$\mathcal{F}$, and that $\mathrm{Im}\,\iota =\mathrm{ker}\,\rho$ in~$\mathcal{F}$. Then, since $\rho$ is a morphism of Lie algebras, we have
\begin{gather*}
\rho [\![\iota (v_1),\iota (v_2)]\!]=[\rho\iota (v_1),\rho\iota (v_2)]=0.
\end{gather*}
Thus, the bracket $[\cdot ,\cdot ]_V$ is well-def\/ined. By a similar argument it can be proved that $\nabla$ and $\mathcal{B}$ are well-def\/ined. The proof of the remaining properties follows the same guidelines that those in Theorem~\ref{tatr}.
\end{proof}

\begin{theorem}\label{iso}
Let
\begin{gather*}
\begin{diagram}
\node{V} \arrow{e,t,J}{\iota} \node{\mathcal{F}} \arrow{e,t,A}{\rho} \node{\mathrm{Der}(\mathcal{A})}
\end{diagram}
\end{gather*}
be a transitive Lie algebroid and $\gamma:\mathrm{Der}(\mathcal{A})\rightarrow \mathcal{F}$ a Lie algebroid connection on $\mathcal{F}$. Then $(\mathcal{F},[\![\cdot ,\cdot ]\!],\rho)$ is isomorphic to $(\mathrm{Der}(\mathcal{A})\oplus V,\langle \cdot ,\cdot \rangle,\mathrm{pr}_1)$, where the bracket is defined through
\begin{gather*}
\langle(X_1,v_1),(X_2,v_2)\rangle =([X_1,X_2],[v_1,v_2]_V+\nabla_{X_1}(v_2)-\nabla_{X_2}(v_1)-\mathcal{B}(X_1,X_2)),
\end{gather*}
and $[\cdot ,\cdot ]_V,\nabla_X$, $\mathcal{B}$ are given in Proposition {\rm \ref{g}}.
\end{theorem}

\begin{proof}
We know that $(\mathrm{Der}(\mathcal{A})\oplus V,\langle \cdot ,\cdot \rangle,\mathrm{pr}_1)$ is a Lie algebroid by Theorem \ref{tat}. Moreover,
\begin{gather*}
\phi: \ \mathrm{Der}(\mathcal{A})\oplus V \longrightarrow  \mathcal{F},\\
\phantom{\phi: \ }{} \ (X,v) \longmapsto  \gamma(X)+\iota(v)
\end{gather*}
is an $\mathcal{A}$-module isomorphism such that
\begin{gather*}
(\rho\circ\phi)(X,v)=\rho(\gamma(X)+\iota(v))=X=\mathrm{pr}_1(X,v)
\end{gather*}
and, by a straightforward computation,
\begin{gather*}
\phi(\langle (X_1,v_1),(X_2,v_2)\rangle )=[\![\phi(X_1,v_1),\phi(X_2,v_2)]\!].\tag*{\qed}
\end{gather*}
\renewcommand{\qed}{}
\end{proof}

\section[Poisson algebras on $\mathcal{A}\oplus V$]{Poisson algebras on $\boldsymbol{\mathcal{A}\oplus V}$}\label{Poisson}

Recall that, if $(\mathcal{A},\{\;,\;\})$ is a Poisson algebra, then we def\/ine for every $f\in\mathcal{A}$ the Hamiltonian derivation $X_f\in \mathrm{Der}(\mathcal{A})$  as the adjoint map $X_f=\{f,\cdot \}$. The following standard properties will be used below:
\begin{enumerate}[(i)]\itemsep=0pt
\item $X_{f_1f_2}=f_1X_{f_2}+f_2X_{f_1}$,
\item the mapping $f\mapsto X_f$ is a Lie algebra morphism, i.e., $X_{\{f_1,f_2\}}=[X_{f_1},X_{f_2}]$.
\end{enumerate}

Now, on $\mathcal{A}\oplus V$ we can def\/ine a product given by
\begin{gather}\label{et2}
(f_1,v_1)\cdot (f_2,v_2)=(f_1f_2,f_1v_2+f_2v_1),
\end{gather}
for $f_1,f_2\in\mathcal{A}$ and $v_1,v_2\in V$.
This makes $\mathcal{A}\oplus V$ a commutative ring. Under certain conditions, this ring is also a Poisson algebra, as shown by the following result.

\begin{theorem}\label{ps}
Let $(\mathcal{A},\{,\}_{\mathcal{A}})$ be a Poisson algebra. Suppose we have a transitive Lie algebroid
\begin{gather*}
\begin{diagram}
\node{V} \arrow{e,t,J}{\iota_2} \node{\mathrm{Der}(\mathcal{A})\oplus V} \arrow{e,t,A}{\rho} \node{\mathrm{Der}(\mathcal{A}),}
\end{diagram}
\end{gather*}
with anchor $\rho=\mathrm{pr}_1$. Then, $\mathcal{A}\oplus V$ is also a Poisson algebra with product defined by \eqref{et2} and Poisson bracket
\begin{gather}\label{cpd}
\{(f_1,v_1),(f_2,v_2)\}=(\{f_1,f_2\}_{\mathcal{A}}, \nabla_{X_{f_1}}(v_2)-\nabla_{X_{f_2}}(v_1)+[v_1,v_2]_V+\mathcal{B}(X_{f_1},X_{f_2})),
\end{gather}
where $[\cdot ,\cdot ]_V$, $\nabla$ and $\mathcal{B}$ are given in Theorem {\rm \ref{tatr}}.
\end{theorem}
\begin{proof}
The Jacobi identity and the Leibniz rule are proved by direct computations, making use of the
$\mathcal{A}$-linearity of the bracket $[\;,\;]_{V}$.
\end{proof}

\begin{corollary}
$V$ is an ideal of $\mathcal{A}\oplus V$ with respect to the Poisson bracket defined by~\eqref{cpd}.
\end{corollary}

The proof of this result follows from the fact $\{0,f\}_{\mathcal{A}}=0$.

\begin{corollary}\label{ps2}
If $(\mathcal{A},\{,\}_{\mathcal{A}})$ is a Poisson algebra and we have a transitive Lie algebroid
\begin{gather*}
\begin{diagram}
\node{V} \arrow{e,t,J}{\iota_2} \node{\mathcal{F}} \arrow{e,t,A}{\rho} \node{\mathrm{Der}(\mathcal{A}),}
\end{diagram}
\end{gather*}
endowed with an algebroid connection $\gamma$, then, $\mathcal{A}\oplus V$ is a Poisson algebra.
\end{corollary}
\begin{proof}
By Theorem \ref{iso}, we know that $(\mathcal{F},[\![\cdot ,\cdot ]\!],\rho)$ is isomorphic to $(\mathrm{Der}(\mathcal{A})\oplus V,\langle\cdot ,\cdot \rangle,\mathrm{pr}_1)$. The statement then follows from Theorem~\ref{ps}.
\end{proof}

\subsection*{Acknowledgements}
JAV and DGB express their gratitude to the members of the Department of Mathematics of the University of Sonora (where part of this work was done), particularly to R.~Flores-Espinoza and G.~D\'avila-Rasc\'on, for their warm hospitality. JAV also thanks the National Council of Science and Technology in Mexico (CONACyT), for the research grant JB2-CB-78791.


\pdfbookmark[1]{References}{ref}
\LastPageEnding

\end{document}